\documentclass[preprint,3p,12pt]{elsarticle}

\usepackage{color}
\usepackage{amssymb}
\usepackage{amsmath}
\usepackage{amsthm}

\usepackage{epsfig}
\usepackage{ulem}
\usepackage{epic}
\usepackage{tabls}
\usepackage{booktabs}
\usepackage{threeparttable}
\usepackage{pict2e}

\usepackage{float}
 \usepackage{subfig}

\usepackage{hyperref}

\hypersetup{backref,colorlinks=true,linkcolor=blue}

\usepackage[vlined]{algorithm2e}

\begin{document}
\begin{frontmatter}



\title{Implicit Methods with Reduced Memory  for Thermal   Radiative Transfer}

\author[ncsu1]{Dmitriy Y. Anistratov}
\author[ncsu2]{Joseph M. Coale}
\address[ncsu]{Department of Nuclear Engineering,
North Carolina State University Raleigh, NC}
\address[ncsu1]{anistratov@ncsu.edu}
\address[ncsu2]{jmcoale@ncsu.edu}

\begin{abstract}
This paper presents  approximation methods for time-dependent thermal
radiative transfer problems in high energy density physics.  It is based on the multilevel quasidiffusion  method defined by the high-order
radiative transfer equation (RTE) and the low-order quasidiffusion  (aka VEF) equations for the moments of the specific intensity.
A large  part of data  storage   in TRT problems between time steps
is determined by  the dimensionality of    grid functions of the radiation intensity.
The approximate implicit methods with reduced memory for  the time-dependent Boltzmann  equation
are applied to the high-order RTE, discretized in time with the backward Euler (BE)  scheme.
The high-dimensional intensity from the previous time level in the BE scheme
 is approximated by means of the low-rank proper orthogonal decomposition (POD).
 Another version of the presented method applies the POD to the remainder term of $P_2$ expansion of the  intensity.
 The accuracy of the solution of the approximate implicit methods depends of the rank of the POD.
  The proposed methods  enable one to reduce   storage requirements in time dependent problems.
 Numerical results of  a  Fleck-Cummings TRT test problem are presented.
\end{abstract}


\end{frontmatter}

\section{Introduction}

We  consider  the  thermal radiative transfer (TRT) problem in 1D slab geometry
that is defined by the time-dependent  radiative transfer equation  (RTE)
\begin{equation} \label{rte}
\frac{1}{c} \frac{\partial I_g  }{\partial t}(x, \mu, t)
+ \mu \frac{\partial  I_g}{\partial  x} (x, \mu, t)
+ \varkappa_{g}(T)I_g(x, \mu, t) =
\varkappa_{g}(T) B_g(T) \, ,
\end{equation}
\[
x \in [0,X] \,  , \quad
\mu \in [-1,1] \, , \quad
 g \in  \mathbb{N}(G) \, , \quad
  t \ge t_0 \, ,
\]
\begin{equation}
\left. I_g \right|_{\stackrel{\mu>0}{x=0}} =  I_g^{in+} \, , \quad
\left. I_g \right|_{\stackrel{\mu<0}{x=X}} =  I_g^{in-} \, , \quad
\left. I_g \right|_{t=t_0} =  I_g^0 \, ,
\end{equation}
and the  material energy balance (MEB) equation
\begin{equation}\label{meb}
\frac{\partial \varepsilon(T)}{\partial t} =
\sum_{g=1}^{G}  \varkappa_{g}(T) \Big(\int_{-1}^1 I_g(x, \mu, t)d  \mu  -   2  B_g(T) \Big)   \, ,
\quad \left. T \right|_{t=t_0} =  T_0 \, ,
\end{equation}
where
$I_g$ is the group specific photon intensity;
$x$ is the spatial position;
$\mu$ is   the direction cosine of particle motion;
$g$ is the index of  photon frequency group;
$\mathbb{N}(G) = \{1,\ldots,G \}$;
$t$ is time;
$\varkappa_g$ is the group opacity;
$T$ is the material temperature;
$\varepsilon$  is the material energy density;
 $B_g$ is the group Planck black-body distribution  function.

The solution of the multigroup RTE in general geometry depends on 7 independent variables.
Temporal discretization schemes for the RTE  involve the discrete solution at the previous time level.
This requires   storing in memory  6-dimensional  grid functions that approximate
the  transport solution on a given mesh in the phase space.
There are  different  approaches  for
 developing   approximate  methods for time-dependent transport problems
that
reduce memory requirements
  \cite{vya-danilova-bnv-1969,matsekh-2018,ryan-mc2019,pg-dya-jctt,dya-tans-2020}.
The $\alpha$-approximation of the intensity in time  reduces   the RTE to
a transport equation of steady-state form with a modified opacity \cite{vya-danilova-bnv-1969}. This approximation assumes
 that the   intensity  varies  exponentially over  each  time interval.
  The approximate
 rate of change in time  can be obtained by means of the solution of low-order moment equations.
 As such,  the $\alpha$-approximation
 rids one of the need to store  the high-dimensional solution from the previous time level
 \cite{vya-danilova-bnv-1969}.
This approximation method for the time-dependent RTE demonstrated good accuracy
 in TRT problems \cite{dya-aristova-vya-mm1996,dya-jcp-2019}.
 Analysis  showed that there are some limitations for  the RTE in the $\alpha$-approximation \cite{pg-dya-jctt}.

Recently,  approximate implicit methods with reduced memory  for  the time-dependent Boltzmann  transport equation have been proposed \cite{dya-tans-2020}.
They use the modified backward Euler (MBE) time integration scheme that applies
 the proper orthogonal decomposition (POD)
 of the  transport solution from the previous time step to compress the data and
 reduce memory requirements  \cite{sirovich-1987,berkooz-holmes-lumley-2005,volkwein-2002}.
The accuracy of the method depends on  the order
  of the low-rank POD of the discrete transport solution.
The error decreases
 as rank increases.  In this paper, we apply the MBE scheme
  within  the  framework of the multilevel quasidiffusion  (MLQD)  method  for solving TRT problems
  \cite{dya-aristova-vya-mm1996,gol'din-1964,gol'din-1972}.

The reminder of the paper is organized as follows.
 In Sec. \ref{sec:method}, the  MLQD  method with approximate implicit scheme is formulated.
 In Sec. \ref{sec:pod}, we present different approximations of the specific intensity by means of the POD.
 The numerical results are presented in Sec. \ref{sec:res}.
 We conclude with a discussion in Sec. \ref{sec:conc}.

\section{\label{sec:method} The MLQD Method  with
	Approximate Implicit  Scheme for the High-Order Problem}

\subsection{MLQD Equations and Discretization}

The  MLQD method is defined by a system of equations consisting of
\begin{enumerate}
\item the multigroup high-order RTE (Eq. \eqref{rte})
\begin{equation} \label{rte-a}
\frac{1}{c} \frac{\partial I_g  }{\partial t}
+ \mu \frac{\partial  I_g}{\partial  x}
+ \varkappa_{g}I_g =
\varkappa_{g} B_g \, ,
\end{equation}
\item the multigroup low-order quasidiffusion (aka VEF) equations for the group radiation energy  density and flux \cite{gol'din-1964,auer-mihalas-1970}
\begin{subequations}\label{mloqd-a}
\begin{equation}
 \frac{\partial E_g   }{\partial t}
+ \frac{\partial F_g   }{\partial x} +
 c\varkappa_{g} E_g    =   2\varkappa_{g} B_g  \, ,
  \end{equation}
  \begin{equation}
 \frac{1}{c}\frac{\partial F_g   }{\partial t} +
  c  \frac{\partial ( f_g    E_g ) }{\partial x}  + \varkappa_{g}  F_g
   =   0   \, ,
\end{equation}
where
\begin{equation}
f_g  = \frac{\int_{-1}^{1} \mu^2 I_g d \mu}{\int_{-1}^{1} I_g d \mu}
\end{equation}
\end{subequations}
is the group QD (Eddington) factor,
\item the effective grey low-order quasidiffusion (LOQD)  equations for the total radiation energy  density and fluxes
 \begin{subequations}\label{gloqd-a}
\begin{equation}\label{gloqd-a-1}
 \frac{\partial E   }{\partial t}
+ \frac{\partial F   }{\partial x}
+ c \bar \varkappa_E  E
= c  \bar \varkappa_B a_R T^4 \, ,
\end{equation}
\begin{equation}\label{gloqd-a-2}
 \frac{1}{c}\frac{\partial F  }{\partial t}
+ c\frac{\partial (\bar f_E  E)   }{\partial x}
+ \bar \varkappa_{| \! F \!|}  F + \bar \eta E = 0 \, ,
\end{equation}
\end{subequations}
 where the spectrum averaged opacities and factors are defined by
 \begin{equation}
 \bar{\alpha}_{H}= \frac{\sum_{g=1}^{G} \alpha_{g} H_g }
{\sum_{g=1}^{G} H_g} \, , \quad
\bar \eta  = \frac{\sum_{g=1}^{G} (\varkappa_{g} - \bar{\varkappa}_{|F|}) F_g}{\sum_{g=1}^{G} E_g} \, ,
\end{equation}
\item the MEB equation \eqref{meb} in  grey form
 \begin{equation}\label{meb-grey-a}
  \frac{\partial \varepsilon(T)  }{\partial t}
 =  c  \big(\bar \varkappa_E E- \bar\varkappa_B a_R T^4 \big) \, .
\end{equation}
\end{enumerate}

 We discretize  the equations of the MLQD method  by the backward Euler (BE) time integration scheme.
  This yields   the semi-discrete RTE  at the $n$-th time level   given by
\begin{equation} \label{rte-be}
\frac{1}{c\Delta t^n} \big(I_{g}^n -   I_{g}^{n-1}\big)  +
\mu \frac{\partial  I_g^n}{\partial  x}
+ \varkappa_{g}^nI_g^n =  Q_g^n\, ,
\end{equation}
where  $\Delta t^n=t^n - t^{n-1}$ is the $n$-th time step, $Q_g^n =\varkappa_{g}(T^n) B_g(T^n)$.
The high-order equation \eqref{rte-be} is discretized  in space by the step characteristic  (SC) scheme.
 The multigroup LOQD equations   discretized in time by the BE   scheme have the following form:
 \begin{subequations}\label{mloqd}
\begin{equation}
\frac{1}{\Delta t^n} \big(E_{g}^n -E_{g}^{n-1}\big)
+ \frac{\partial F_g^n   }{\partial x} +
 c\varkappa_{g}^n E_g^n    =   2Q_g^n  \, ,
  \end{equation}
  \begin{equation}
\frac{1}{c\Delta t^n} \big(F_{g}^n - F_{g}^{n-1}\big)  +
  c  \frac{\partial ( f_g^n    E_g^n ) }{\partial x}  + \varkappa_{g}^n  F_g^n
   =   0   \, ,
\end{equation}
\end{subequations}
\begin{equation} \label{f}
f_g^n  = \frac{\int_{-1}^1 \mu^2   I_g^n d \mu}{\int_{-1}^1 I_g^n d \mu}  \, .
\end{equation}
 The grey LOQD and MEB equations approximated with the BE scheme are defined by
 \begin{subequations}\label{gloqd}
\begin{equation}\label{gloqd-1}
\frac{1}{\Delta t^n} \big(E^n -E^{n-1}\big)
+ \frac{\partial F^n   }{\partial x}
+ c \bar \varkappa_E^n  E^n
= c  \bar \varkappa_B^n a_R (T^n)^4 \, ,
\end{equation}
\begin{equation}\label{gloqd-2}
\frac{1}{c\Delta t^n} \big(F^n -F^{n-1}\big)
+ c\frac{\partial (\bar f_E^n  E^n)   }{\partial x}
+ \bar \varkappa_{| \! F \! |}^n  F^n + \bar \eta^n E^n = 0 \, .
\end{equation}
\end{subequations}
\begin{equation}\label{meb-grey}
\frac{1}{\Delta t^n} \big(\varepsilon(T^n)-\varepsilon(T^{n-1})\big)
 =  c  \big(\bar \varkappa_E^n E^n - \bar\varkappa_B^n a_R (T^n)^4\big) \, .
\end{equation}
  The multigroup LOQD  equations  are discretized in space by  a second-order finite volume (FV) method. The spatial discretization of the grey LOQD equations   is algebraically  consistent with  the  discretized multigroup LOQD equations \cite{dya-jcp-2019}.
  We refer to the described method as the MLQD method with BE-SC scheme.

\subsection{Approximate Implicit Method for the RTE}

In the approximate implicit scheme,
 the multigroup RTE    \eqref{rte-a} is discretized   by the MBE time integration scheme given by  \cite{dya-tans-2020}
\begin{equation} \label{rte-mbe}
\frac{1}{c\Delta t^n} \big(I_{g}^n - \hat I_{g}^{n-1}\big)  +
\mu \frac{\partial  I_g^n}{\partial  x}
+ \varkappa_{g}^nI_g^n =  Q_g^n\, ,
\end{equation}
where the grid functions of group intensity  $\hat I_{g}^{n-1}$  are approximated by
 the low-rank  POD  of  the solution $I_{g}^{n-1}$ computed at  the  time step $n-1$.
The SC scheme for the high-order equation \eqref{rte-mbe}
  is formulated for the cell-edge  $\big(I_{g \, m \, j+1/2}^n \big)$  and cell-average $\big( I_{g \, m \, j}^n \big)$  angular fluxes by  means of the detailed particle
  balance equation and weighted auxiliary relation
  \begin{subequations} \label{mbe-sc}
  \begin{equation} \label{sc-bal}
{\frac{\Delta x_j}{c\Delta t^n}}
\big(I_{g \, m \, j}^n  - \hat I_{g \, m \, j}^{n-1}\big)  +\mu_m \big( I_{g \, m \, j+1/2}^n  -  I_{g \, m \, j-1/2}^n \big)
+ \varkappa_{g \, j}^n I_{g \, m \, j}^n  \Delta x_j  =
Q_{g  \, j}^n \Delta x_j   \, ,
\end{equation}
\begin{equation} \label{sc-aux}
I_{g \, m \, j}^n =  \gamma_{g \, m \, j}^n   I_{g \, m \, j-1/2}^n   + (1 - \gamma_{g \, m \, j}^n) I_{g \, m \, j+1/2}^n \, ,
\end{equation}
\begin{equation}
 \gamma_{g \, m \, j}^n = \frac{1}{\tau_{g \, m \, j}^n} - \frac{1}{e^{\tau_{g \, m \, j}^n}  - 1} \, , \quad
 \tau_{g \, m \, j}^n = \frac{1}{\mu_m}\Big(\varkappa_{g \, j}^n + (c\Delta t^n)^{-1} \Big) \Delta x_j  \, ,
\end{equation}
\end{subequations}
where $m \in \mathbb{N}(M)$ is the index of angular direction,
 $j  \in \mathbb{N}(J)$ is the index of the spatial interval, $\Delta x_j$ is the width of the $j$-th cell.
We refer to the discretized   RTE  \eqref{mbe-sc} as the MBE-SC scheme.

\section{\label{sec:pod}  Approximation of the Specific Intensity}

\subsection{POD of the Intensity}

The  MBE-SC scheme (Eqs. \eqref{mbe-sc}) needs to store the cell-average intensity $I_{g \, m \, j}^n$.
In each photon frequency group, it is   a 2D discrete grid function of $j$ and $m$.
We interpret it  in a  matrix  form  defined by
$A_{I} = [\boldsymbol{I}_{1} \ldots \boldsymbol{I}_{M}]$
  $\big(A_{I} \in \mathbb{R}^{J \times M} \big)$, where
the columns    are given by
$\boldsymbol{I}_{m}  = [I_{\, m \, 1} \ldots I_{m \, J}]^T$
 $\big(\boldsymbol{I}_{m}  \in \mathbb{R}^J\big)$.
  Here we omitted group and time indices for the sake of brevity.
We approximate  the grid function of the group cell-average intensity by the low--rank POD
\cite{berkooz-holmes-lumley-2005,volkwein-2002}.
 The reduced  singular value decomposition (SVD) of $A_{I} $  has the form:
\begin{equation} \label{A-svd}
A_{I}  =  U_{I}    \Lambda_{I}  V_{I} ^T \, .
\end{equation}
$\Lambda_{I}  =  \text{diag}(\lambda_1 \ldots \lambda_d) \ \in \mathbb{R}^{d \times d }$
is  the diagonal matrix of singular values,
where
\begin{equation} \label{d}
d= \min(J,M)
\end{equation}
 is the rank of $A_{I} $.
$U_{I}  =  [ \mathbf{u}_1 \ldots \mathbf{u}_d]$
and $V_{I}   =  [ \mathbf{v}_1 \ldots \mathbf{v}_d]$
are the matrices  of left and right singular vectors, respectively,
where
$\mathbf{u}_{\ell}   \in \mathbb{R}^J$ and $\mathbf{v}_{\ell}    \in \mathbb{R}^M$.
The approximate group intensity
\linebreak
$\boldsymbol{ \hat  I}_{m} =[\hat I_{ m \, 1}  \ldots  \hat I_{ m \, J}]^T$ is defined by the low-rank POD   of $A_{I}$ given by
\begin{equation} \label{POD-I}
\hat A_{I}^{r}    = \sum_{\ell = 1}^r \lambda_{\ell} \mathbf{u}_{\ell}  \otimes (\mathbf{v}_{\ell})^T \, , \quad  r < d \, , \quad
\mbox{where} \quad
  \hat A_{I}^{r} = [\boldsymbol{\hat I}_{1}  \ldots  \boldsymbol{\hat I}_{M}] \, .
\end{equation}
This is the optimal approximation of the matrix  $A_{I}$ in the 2-norm \cite{volkwein-2002,ipsen-2009}.
The low-rank approximation  \eqref{POD-I} requires  storage of
the first $r$ singular values and associated  left and right singular vectors.
Thus, this  approximation leads to  memory allocation of a data set with the number of elements
\linebreak
 $r(J+M+1)$  in each group.
The rank can be chosen according to various criteria.

\subsection{POD of the Remainder Term}

We cast the intensity  as its  $P_2$ approximation and the remainder term
defined by
\begin{equation}
\Delta I_{m \,j} = I_{m \,j} -
\frac{1}{2} \Big(  \tilde \phi_{j} + 3 \mu_m \tilde F_{j}
+ \frac{5}{4}  \big( 3 \mu_m^2 - 1 \big) \big(3   f_{j} - 1 \big) \tilde \phi_{j}  \Big)  \, ,
\end{equation}
where the $P_2$ expansion coefficients are
calculated by the solution of the high-order RTE, namely,
\begin{equation}
\mathcal{ \tilde \phi}_{j} = \sum_{m=1}^M I_{m \,j} w_m \, , \quad
  \tilde F_{j} = \sum_{m=1}^M \mu_m I_{m \,j} w_m \, , \quad
f_{j} = \frac{\sum_{m=1}^M \mu_m^2 I_{m \,j} w_m}
{\sum_{m=1}^M   I_{m \,j}^n w_m} \, .
\end{equation}
The discrete 2D function $\Delta I_{m \,j}$ is treated as  a matrix defined by
$\Delta  A_{I} = [\Delta \boldsymbol{I}_{1}  \ldots  \Delta\boldsymbol{I}_{M}  ]$,
where
$\Delta  \boldsymbol{I}_{m} = [\Delta I_{ m \, 1}  \ldots \Delta I_{m \, J}]^T$.
Its POD is given by
\begin{equation} \label{delta-A-svd}
\Delta A_{I}   =  U_{I}^{\prime}   \Lambda_{I}^{\prime} (V_{I}^{\prime} ) ^T \, .
\end{equation}
where
$\Lambda_{I}^{\prime}  =  \text{diag}(\lambda_1^{\prime} \ldots  \lambda_d^{\prime})\ \in \mathbb{R}^{d \times d }$,
$U_{I}^{\prime}  =  [ \mathbf{u}_1^{\prime}  \ldots \mathbf{u}_d^{\prime}]$,
$V_{I}^{\prime}   =  [ \mathbf{v}_1^{\prime} \ldots \mathbf{v}_d^{\prime}]$,
$\mathbf{u}_{\ell}^{\prime}   \in \mathbb{R}^J$,  and $\mathbf{v}_{\ell}^{\prime} \in \mathbb{R}^M$.
We apply the low-rank POD
\begin{equation}
\Delta \hat  A_{I}^r   = [\Delta \boldsymbol{\hat I}_{1}  \ldots  \Delta \boldsymbol{\hat I}_{M} ]  = \sum_{k = 1}^r \lambda_{k}^{\prime} \mathbf{u}_{k}^{\prime}  \otimes (\mathbf{v}_{k}^{\prime})^T \, , \  \
  r < d \, ,
\ \
 \Delta  \boldsymbol{\hat I}_{m} = [\Delta \hat I_{ m \, 1}, \ldots, \Delta \hat I_{m \, J}]^T
 \end{equation}
to define approximate intensities as the sum of  its  $P_2$ approximation and the POD of the remainder term
\begin{equation} \label{POD-delta-I}
\hat I_{m \,j} =
\frac{1}{2} \Big(  \tilde \phi_{j} + 3 \mu_m \tilde F_{j}
+ \frac{5}{4}  \big( 3 \mu_m^2 - 1 \big) \big(3   f_{j}^n - 1 \big) \tilde \phi_{j}  \Big)
+   \Delta \hat I _{m \,j}\, .
\end{equation}
This  approximation needs to store in memory
 $r(J+M+1) + 2J$ elements that includes (i)  $r(J+M+1) $ elements for the remainder term and (ii) $2J$ elements for vectors of two angular moments $\boldsymbol{\tilde \phi}$ and $\boldsymbol{\tilde F}$.

\section{\label{sec:res} Numerical Results}

We present   numerical results of the   Fleck-Cummings (F-C) test \cite{fleck-1971}.
The spatial domain
($0 \le x\le 6$) contains one material.
The  spectral opacity of the material  is   given by \linebreak
$\varkappa_{\nu} = \frac{27}{(h\nu)^3}\big(1-e^{-\frac{h\nu}{kT}}\big)$.
There is
incoming radiation with black-body spectrum $B_{\nu}$ at temperature
$kT_{in}=1$ keV at the left boundary.
The right boundary is vacuum.
The initial temperature of the slab is  $kT_0=1$ eV.
At $t=0$  the radiation intensity   in the slab   has   the black-body spectrum at $T_0$.
The material energy  density  is  $\varepsilon = c_\nu T$, where  $c_\nu= 0.5917a_RT_{in}^3$.
 The problem is solved over the time interval  $0 \le t  \le 6$ ns.
 The time step  size is  $\Delta t~= 2 \times 10^{-2}$~ns.
   The uniform spatial mesh consists of  $J=100$ cells.
The angular mesh has  8  discrete directions ($M=8$).
The double $S_4$ Gauss-Legendre  quadrature set is used.
We define $G=17$ energy groups.   The parameters of convergence criteria for temperature and energy density are  $\epsilon_T=\epsilon_E=10^{-12}$, respectively.

\begin{figure}[h!]
\centering
\subfloat[$\frac{||T_h - T_h^r||_{\infty}}{||T_h||_{\infty}}$. \label{MLQD-POD-I-T}]{\includegraphics[scale=0.32]{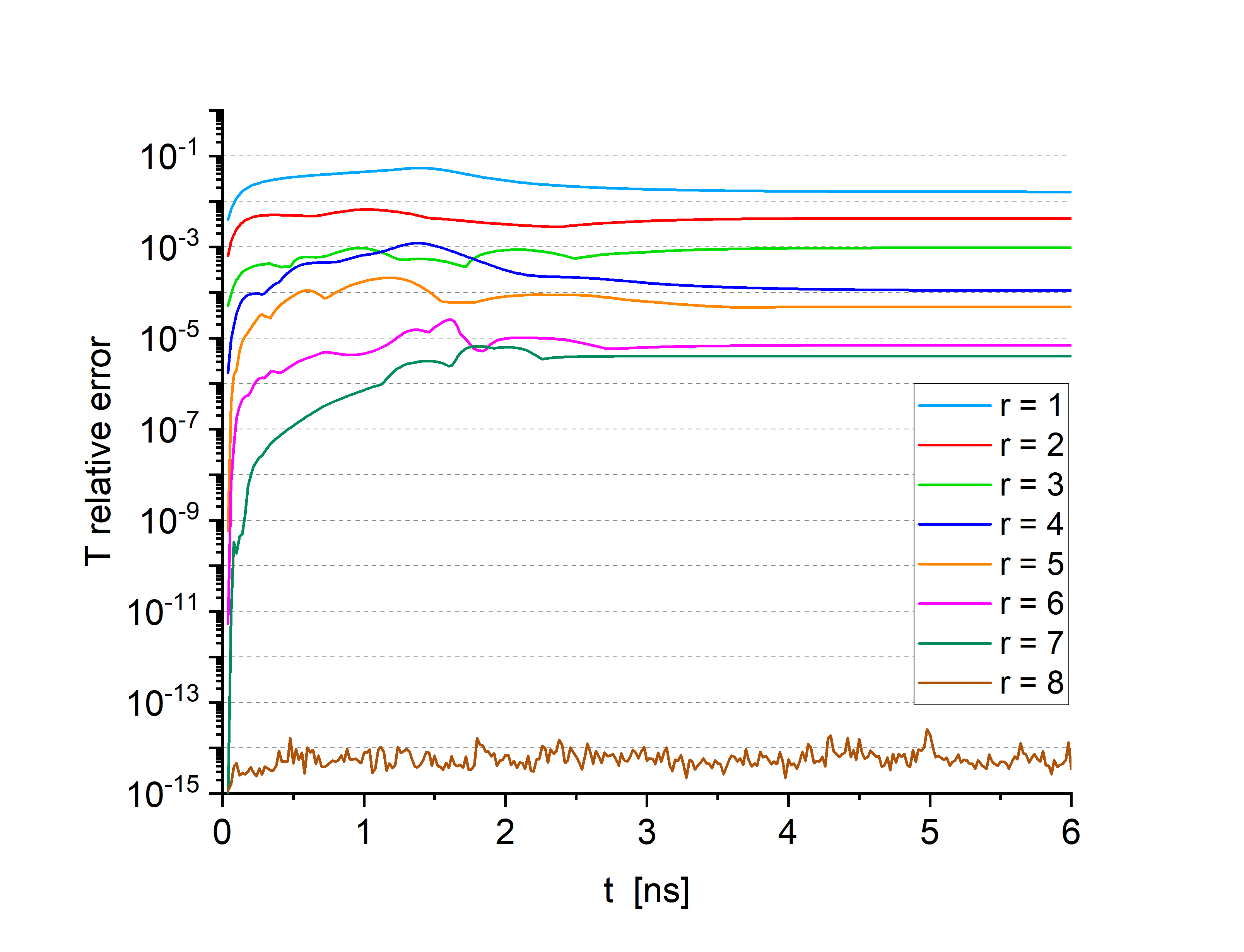}}
\subfloat[$\frac{||E_h - E_h^r||_{\infty}}{||E_h||_{\infty}}$. \label{MLQD-POD-I-E}]{\includegraphics[scale=0.32]{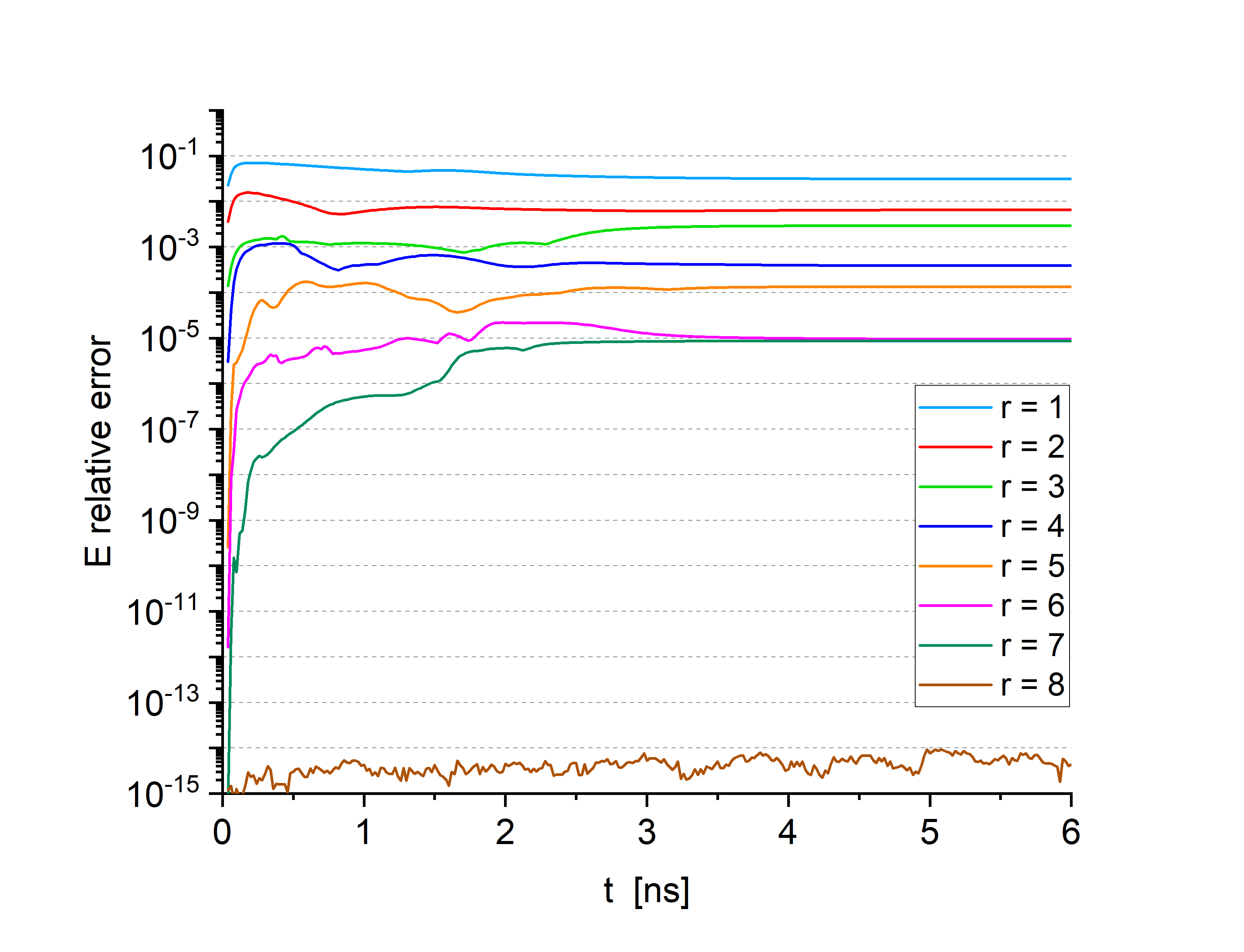}}
\caption{\label{MLQD-POD-I}
Relative error   in $\infty$-norm  of the solution of  the MLQD method  with the MBE-SC scheme and POD of the intensity compared to the discrete solution on the corresponding grid in phase space and time.}

\centering
\subfloat[$\frac{||T_h - T_h^r||_{\infty}}{||T_h||_{\infty}}$. \label{MLQD-POD-delta-I-T}]{\includegraphics[scale=0.32]{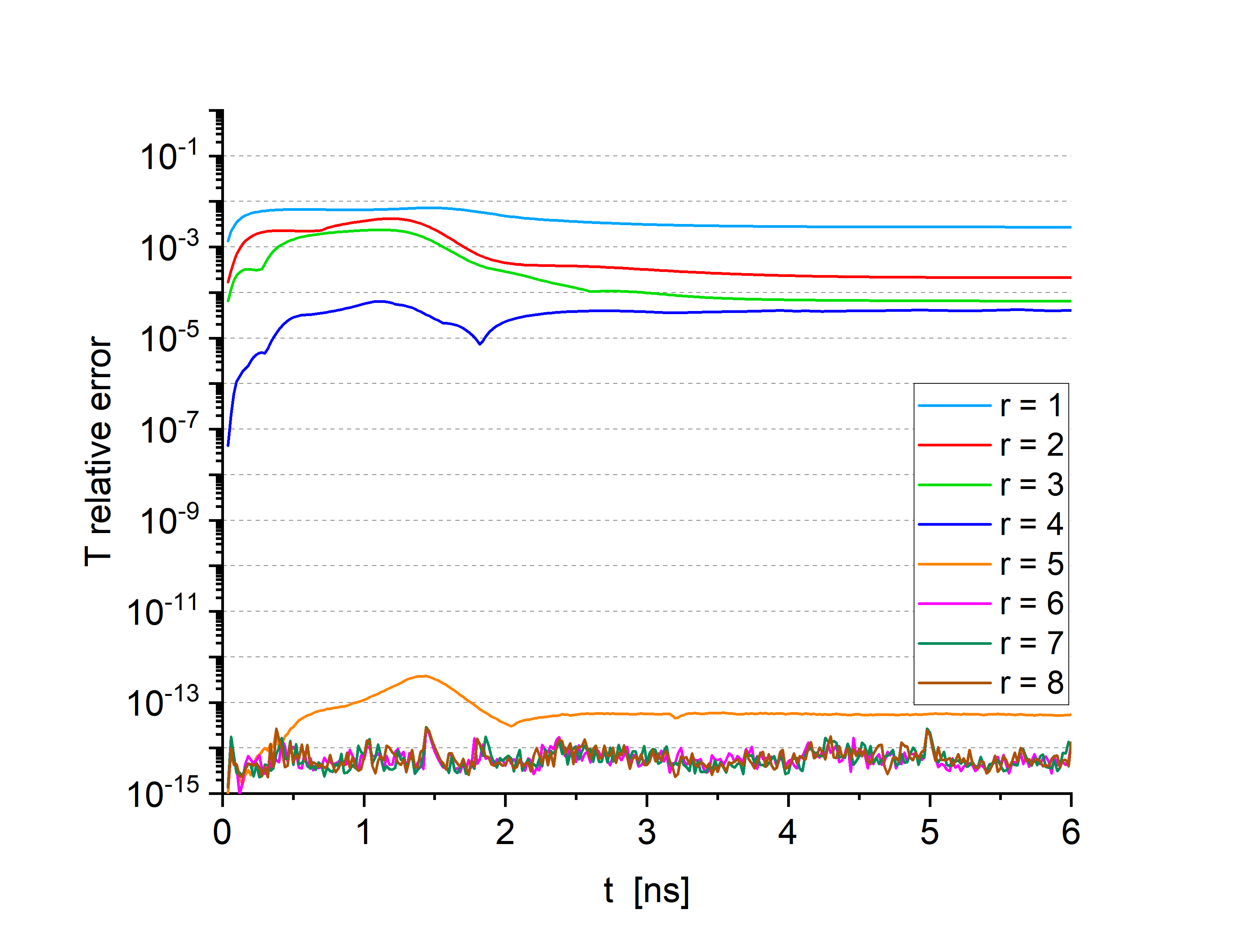}}
\subfloat[$\frac{||E_h - E_h^r||_{\infty}}{||E_h||_{\infty}}$. \label{MLQD-POD-delta-I-E}]{\includegraphics[scale=0.32]{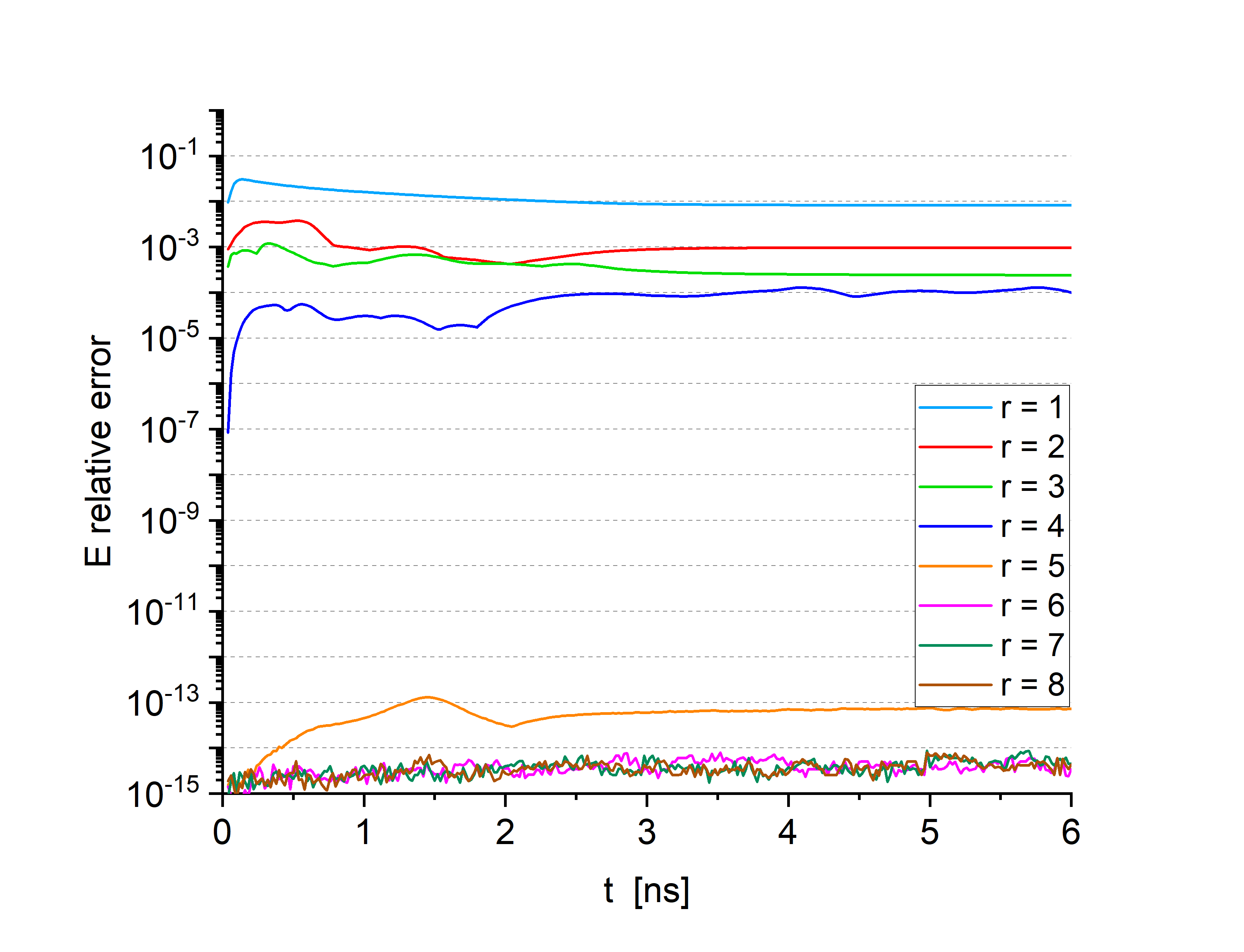}}
\caption{\label{MLQD-POD-delta-I}
Relative error   in $\infty$-norm  of the solution of  the MLQD method with the MBE-SC scheme and POD of the remainder term compared to the discrete solution on the corresponding grid in phase space and time.}
\end{figure}

\begin{figure}[t]
\centering
\subfloat[$\frac{||T_h - T_h^r||_{\infty}^{\text{POD-RT}}}{||T_h - T_h^r||_{\infty}^{\text{POD-I}}}$. \label{ratio-err-T}]{\includegraphics[scale=0.32]{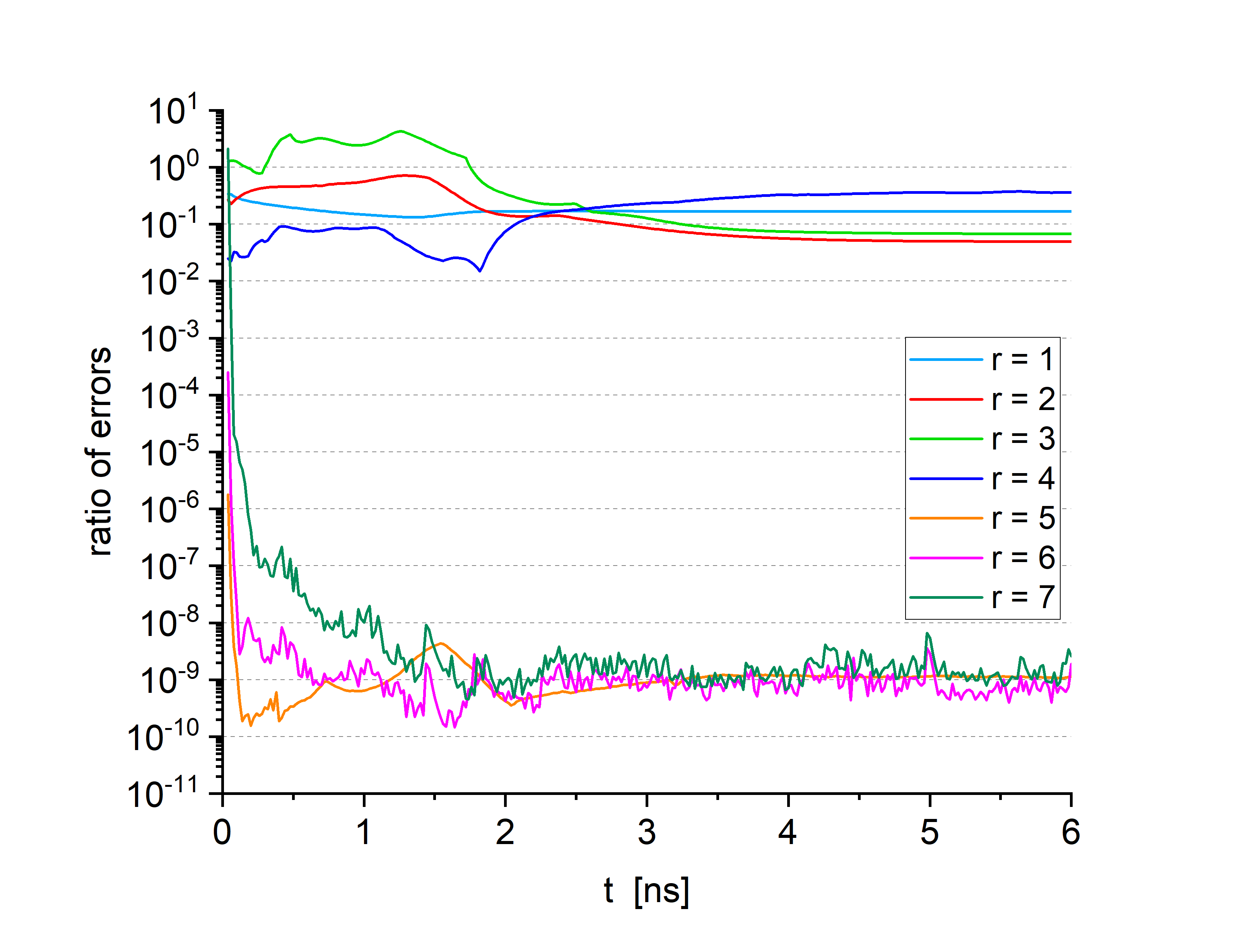}}
\subfloat[$\frac{||E_h - E_h^r||_{\infty}^{\text{POD-RT}}}{||E_h - E_h^r||_{\infty}^{\text{POD-I}}}$. \label{ratio-err-E}]{\includegraphics[scale=0.32]{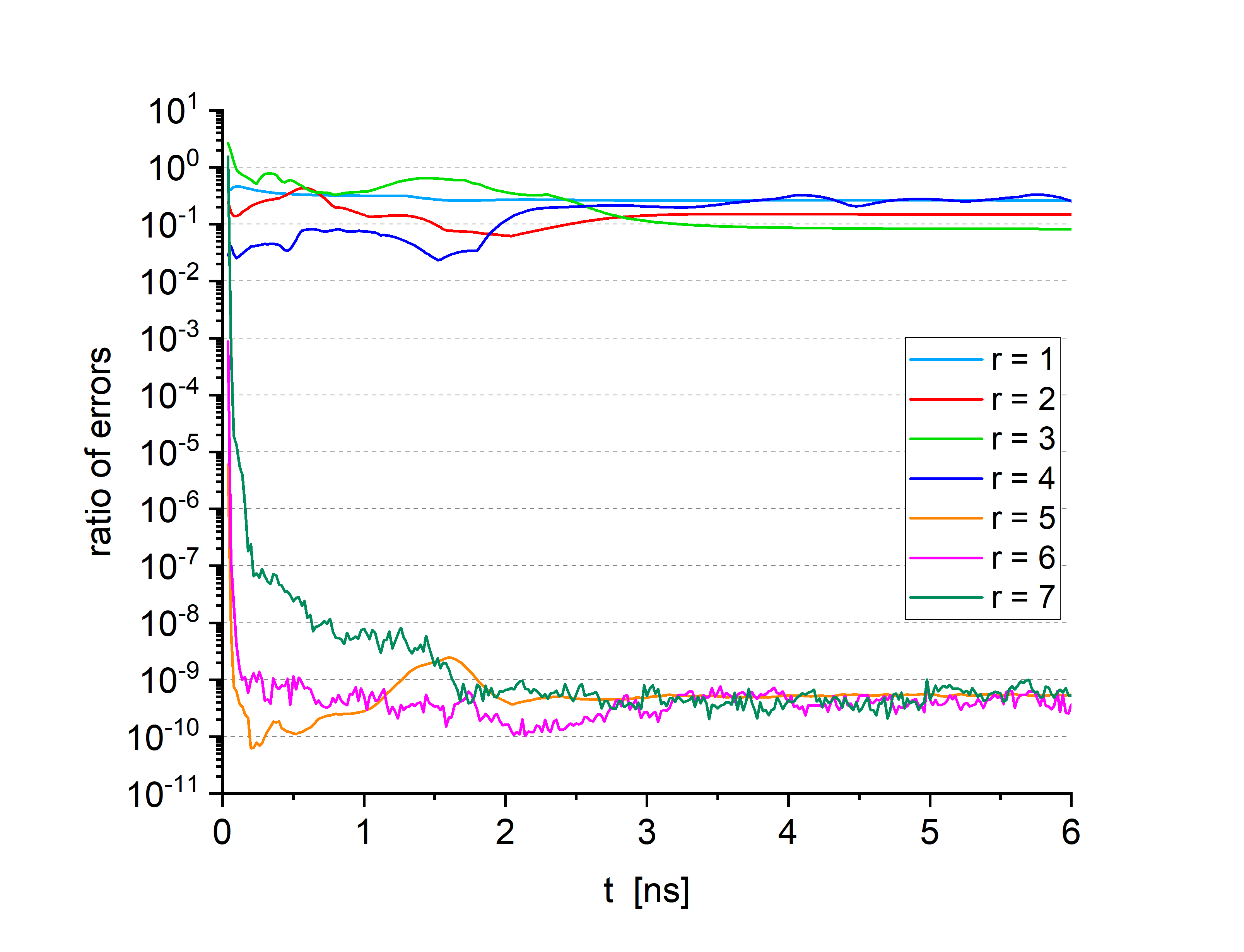}}
\caption{\label{ratio-err}
The  error  of the method with POD  of the remainder term (POD-RT)  over the error of the method with POD of intensity (POD-I). }
\end{figure}

 The discrete solution of the MLQD method with the MBE-SC scheme, namely,
 the total radiation energy density $E^r_h$ and  temperature   $T^r_h$
 of the approximate implicit method with the rank $r$  POD
 is compared to
  the discrete solution $T_h$ and $E_h$  of  the MLQD method with the BE-SC scheme on the corresponding  grid in the phase space and time.
 The numerical results   of  the method with  the  POD of the intensity of the rank $r$ in all groups    are presented in   Figure  \ref{MLQD-POD-I}.
  The plots show the relative error  in reproducing the discrete solution in $\infty$-norm, namely,
   $\frac{||T_h - T_h^r||_{\infty}}{||T_h||_{\infty}}$ and
 $\frac{||E_h - E_h^r||_{\infty}}{||E_h||_{\infty}}$
  for the complete range of $r$.
   The  results   obtained with  the  MBE-SC scheme
  using the  POD of  the rank $r$  of the remainder term  in each group   are shown in   Figure  \ref{MLQD-POD-delta-I}.   In this test, the full rank $d$ (Eq. \eqref{d})  equals 8.
  The results with the full-rank POD ($r=8$) of both methods   illustrate that they accurately reproduce
  the discrete solution of the MLQD method with the SC scheme on the given grid as expected.
In case $r=5,6,7$  the solution of the method with the POD of the remainder term has very small error.
This is due to  explicit  accounting for the first three  Legendre moments of the intensity (Eq. \eqref{POD-delta-I}). The singular eigenvalues $\lambda_{\ell}^{\prime}$ for $\ell=5,6,7$ in groups are very small.
In this test problem, the method with POD of the remainder term is predominantly more accurate  than to the method with POD of the intensity for the given rank $r$.  However, it uses more data for the rank $r$.
Figure \ref{ratio-err} shows the ratio between errors of the method with  the POD of the remainder term (POD-RT)
and the one  with the POD of the intensity (POD-I).

The gains in memory allocation  depend on both the number of spatial cells $J$ and angular directions $M$ and hence are problem specific.
For the phase-space grid used in the test,
   the size of the data set stored by this  MLQD method with
  the RTE discretized the BE-SC  scheme  at the end of each time step is
  $D=G(J \times M + 2 \times J +1 ) + 2\times J +1=17218$.
  This includes the data for (i) the multigroup RTE, (ii) the multigroup and grey LOQD equations,
   and (iii) the MEB equation.
Table~\ref{data-POD} shows the percentage reduction of required data storage sizes
of  the MLQD method with
each of the two versions of
the  MBE-SC scheme compared to that of the MLQD method with the BE-SC scheme.
 Negative values indicate an increase in storage compared to the BE-SC scheme.
 In this test, the method with POD of intensities shows gains in memory for all ranks, i.e. $r=1,\ldots,7$.  The method with POD of the remainder term reduces memory allocation for
  $r=1,\ldots,5$.

 \begin{table}[h]
 	\centering
	\caption{\label{data-POD} Reduction [\%] in memory storage of previous step data
   of the MLQD method with the MBE-SC scheme ($J$=100, $M$=8).}
   \begin{tabular}{|c|c|c|c|c|c|c|c|}
  \hline
 Rank ($r$)  & 1 & 2 & 3 & 4 & 5 & 6 & 7 \\ \hline
POD-I& 68.2 &  57.5 &  46.7 &  35.9 &  25.2 &  14.4 &  3.7  \\ \hline
POD-RT & 48.5 & 37.7	 & 27.0 & 16.2	 & 5.4 & -5.3 & -16.1  \\ \hline
  \end{tabular}
 \end{table}

\begin{figure}[h!]
	\centering \hspace*{-.5cm}
	\subfloat[$\frac{||E_h - E_h^r||_{\infty}}{||E_h||_{\infty}}$ at $t=0.4$ ns]{\includegraphics[width=.35\textwidth]{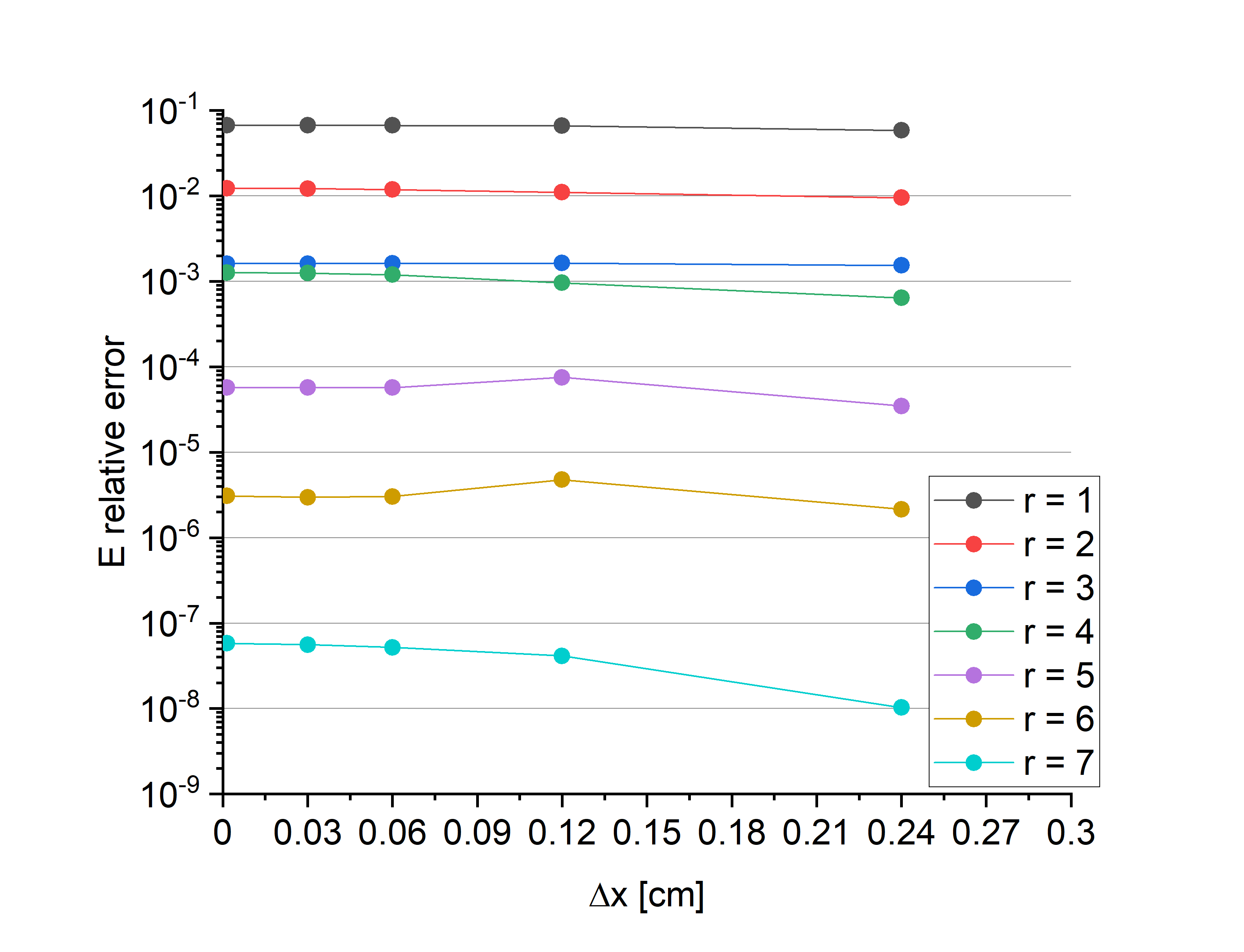}}
	\subfloat[$\frac{||E_h - E_h^r||_{\infty}}{||E_h||_{\infty}}$ at $t=1$ ns]{\includegraphics[width=.35\textwidth]{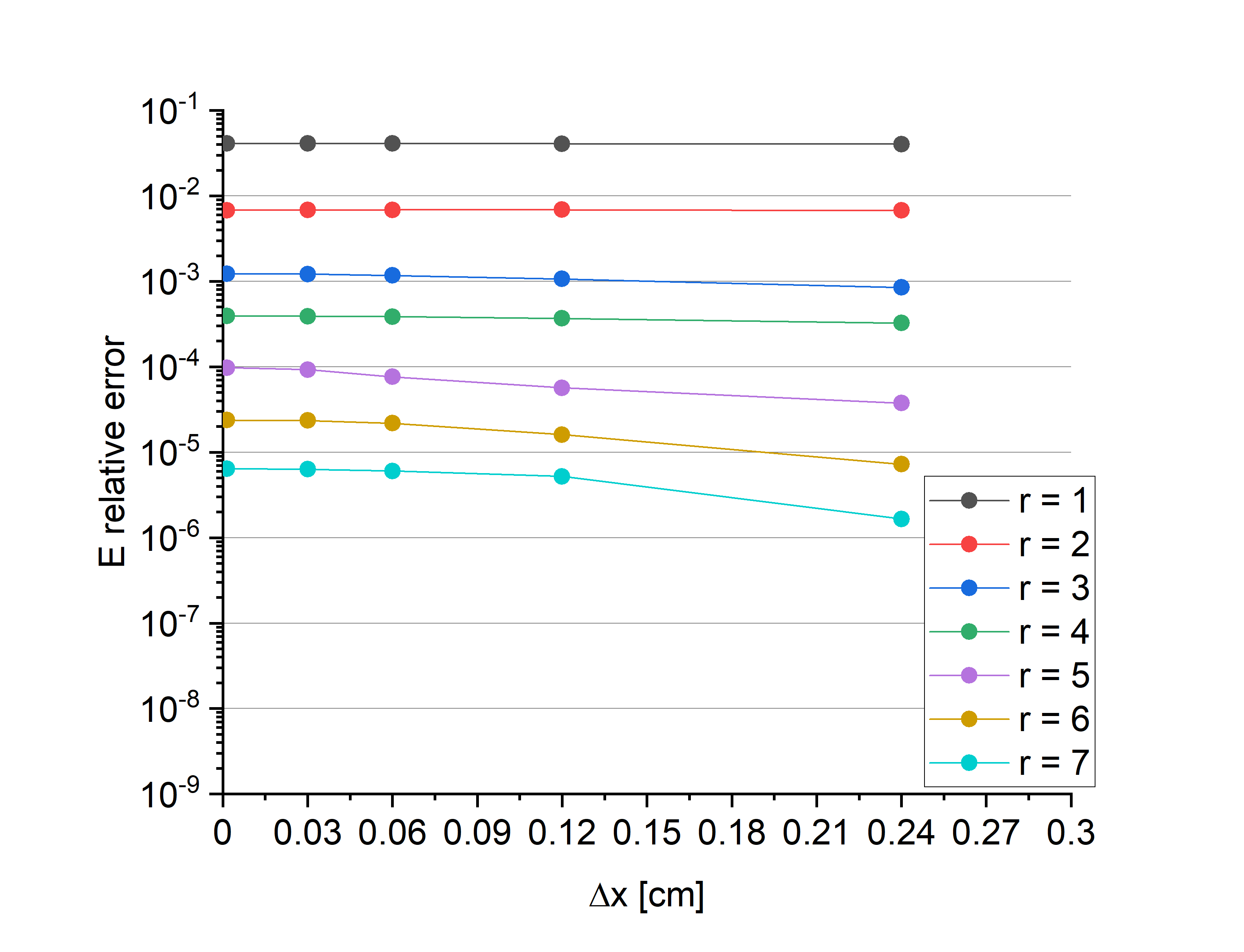}}
	\subfloat[$\frac{||E_h - E_h^r||_{\infty}}{||E_h||_{\infty}}$ at $t=6$ ns ]{\includegraphics[width=.35\textwidth]{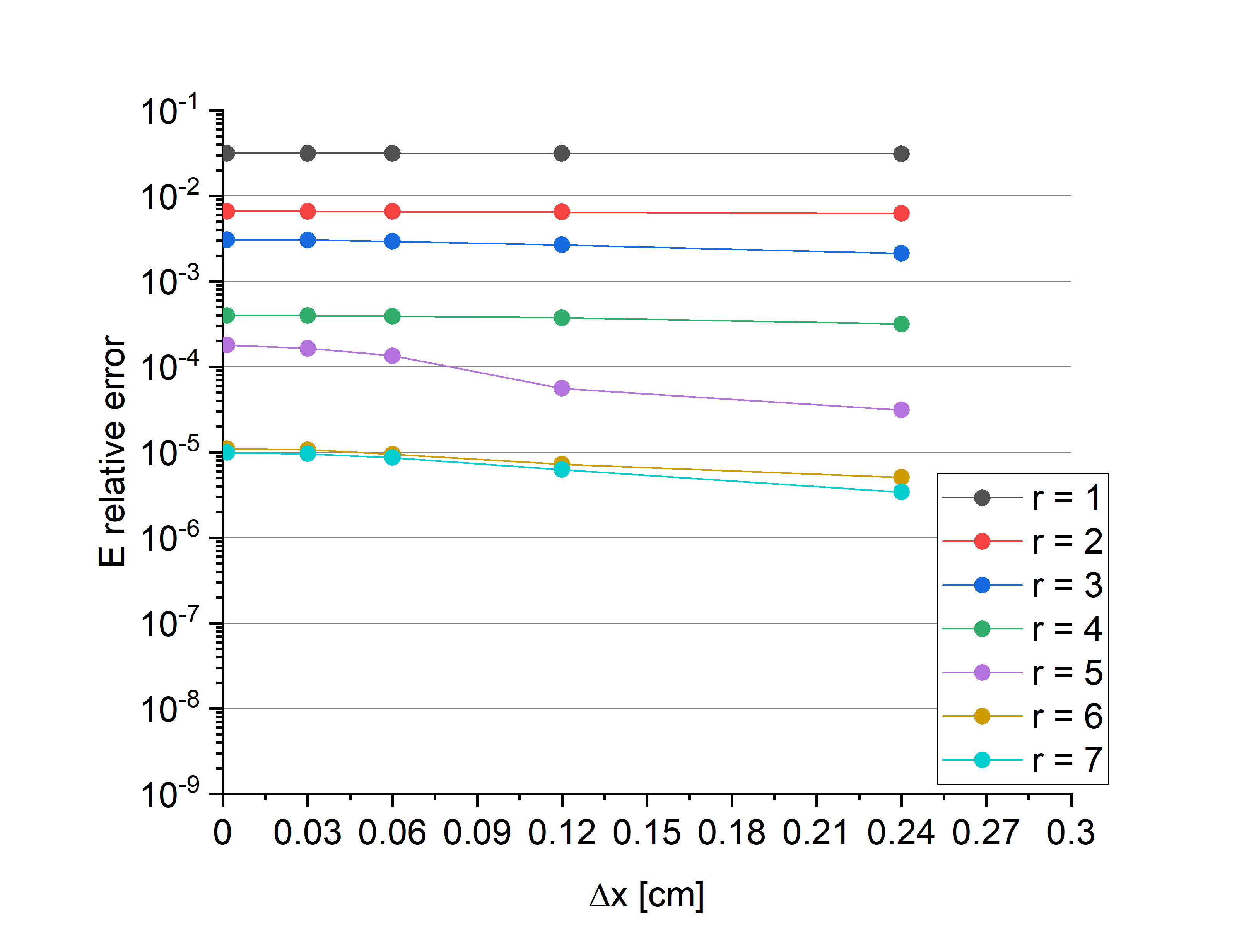}}
	\caption{Results of refinement of spatial mesh for
 the MLQD method with the MBE-SC scheme and POD of the intensity for $\Delta t = 2\times 10^{-2}$.
\label{fig:refine-x}}

	\centering \hspace*{-.5cm}
	\subfloat[$\frac{||E_h - E_h^r||_{\infty}}{||E_h||_{\infty}}$ at $t=0.4$ ns]{\includegraphics[width=.35\textwidth]{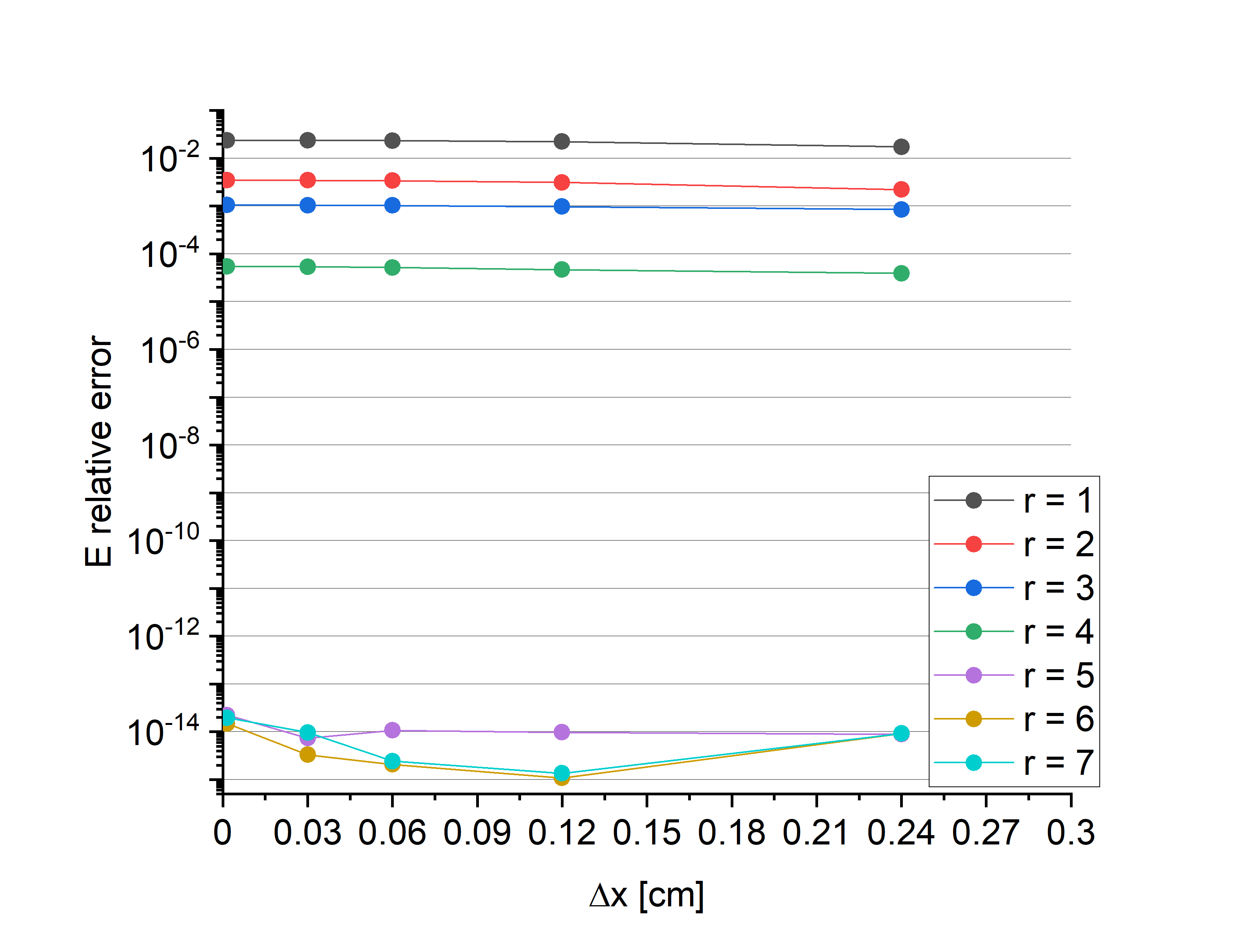}}
	\subfloat[$\frac{||E_h - E_h^r||_{\infty}}{||E_h||_{\infty}}$ at $t=1$ ns]{\includegraphics[width=.35\textwidth]{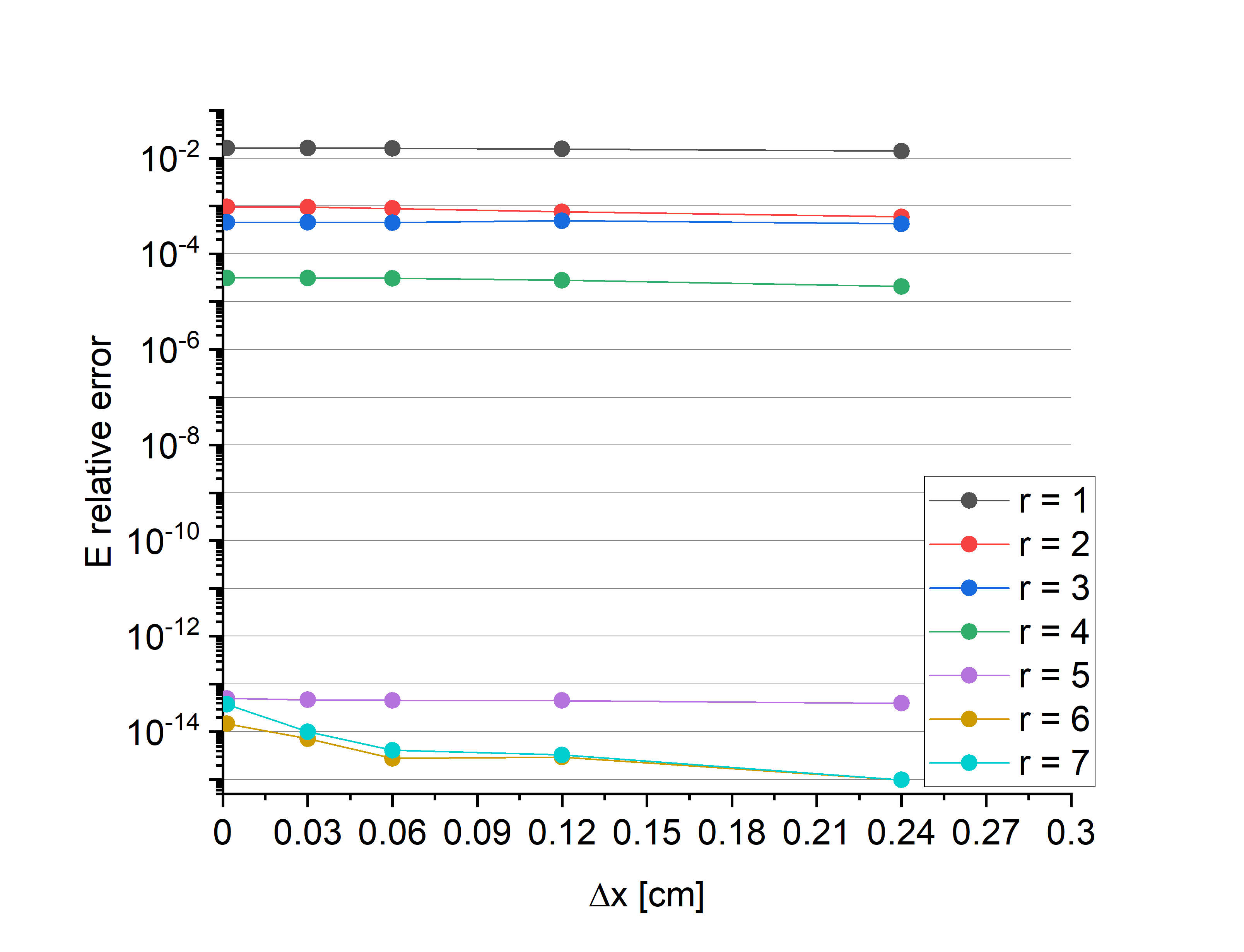}}
	\subfloat[$\frac{||E_h - E_h^r||_{\infty}}{||E_h||_{\infty}}$ at  $t=6$ ns ]{\includegraphics[width=.35\textwidth]{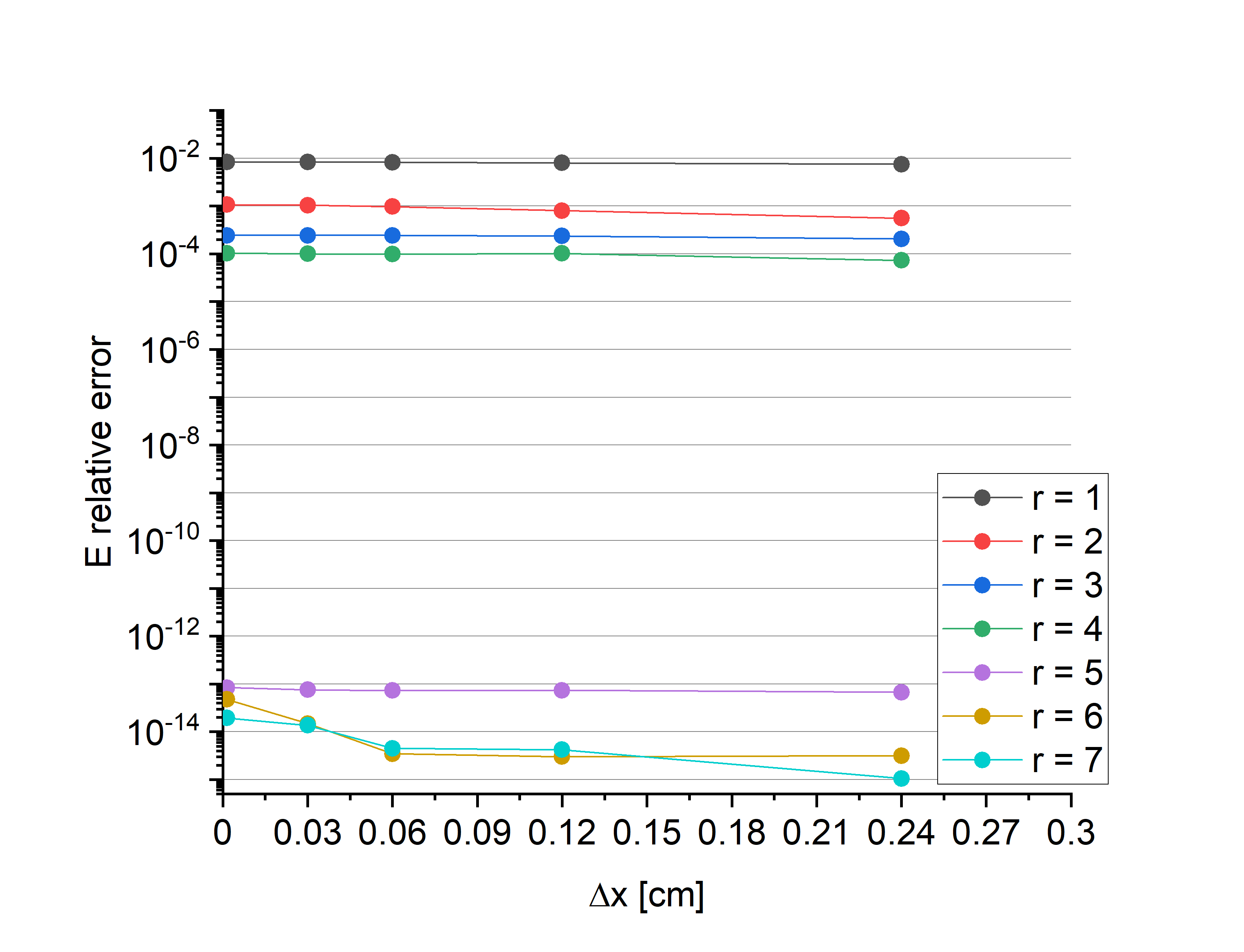}}
	\caption{Results of spatial mesh  refinement   for
the MLQD method with the MBE-SC scheme and POD of the remainder term
for $\Delta t = 2\times 10^{-2}$ ns. \label{fig:refine-x-P2}}
\vspace{-0.5cm}
\end{figure}

\begin{figure}[t]
	\vspace*{-1cm}
	\centering \hspace*{-.5cm}
	\subfloat[$\frac{||E_h - E_h^r||_{\infty}}{||E_h||_{\infty}}$ at $t=0.4$ ns]{\includegraphics[width=.35\textwidth]{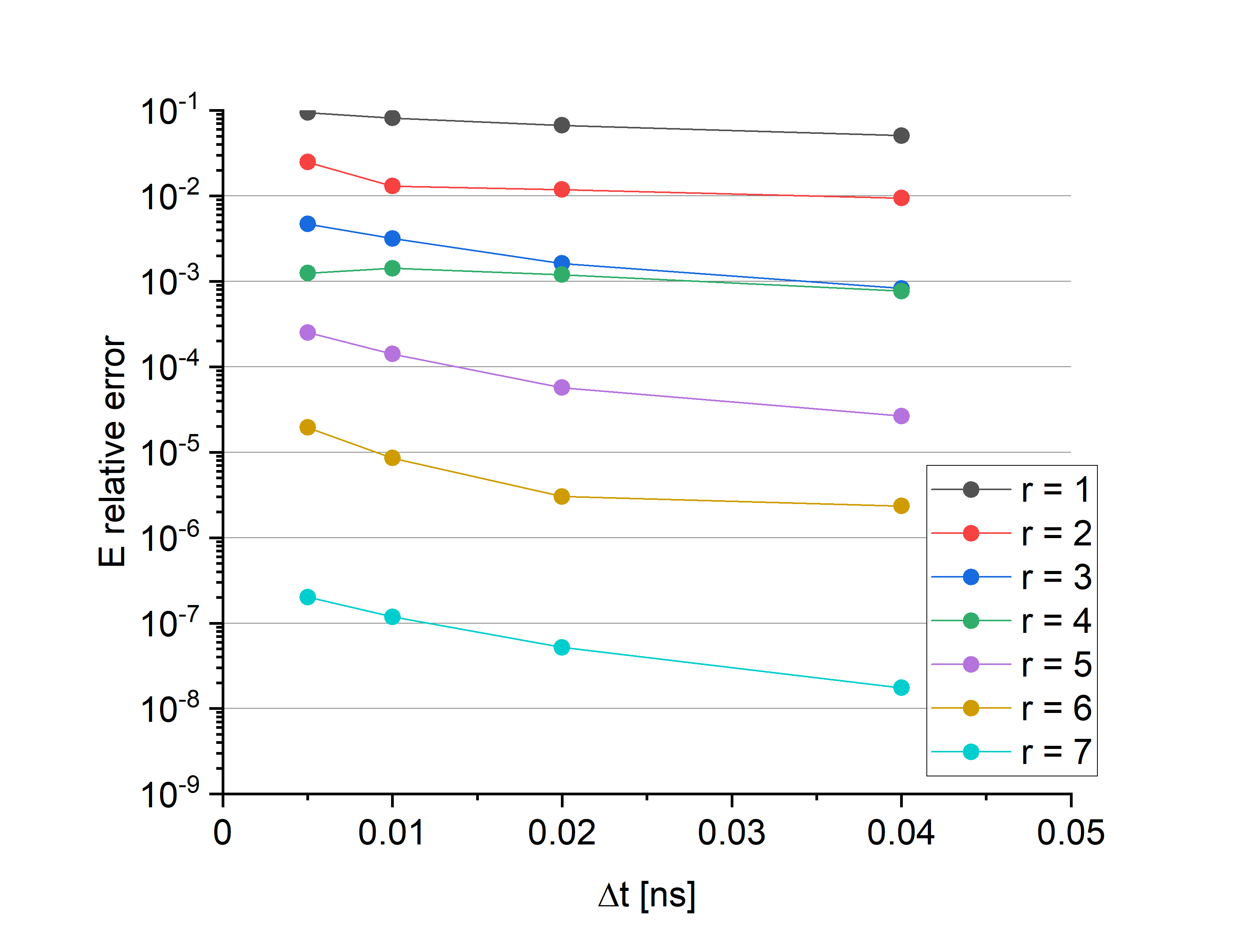}}
	\subfloat[$\frac{||E_h - E_h^r||_{\infty}}{||E_h||_{\infty}}$ at $t=1$ ns]{\includegraphics[width=.35\textwidth]{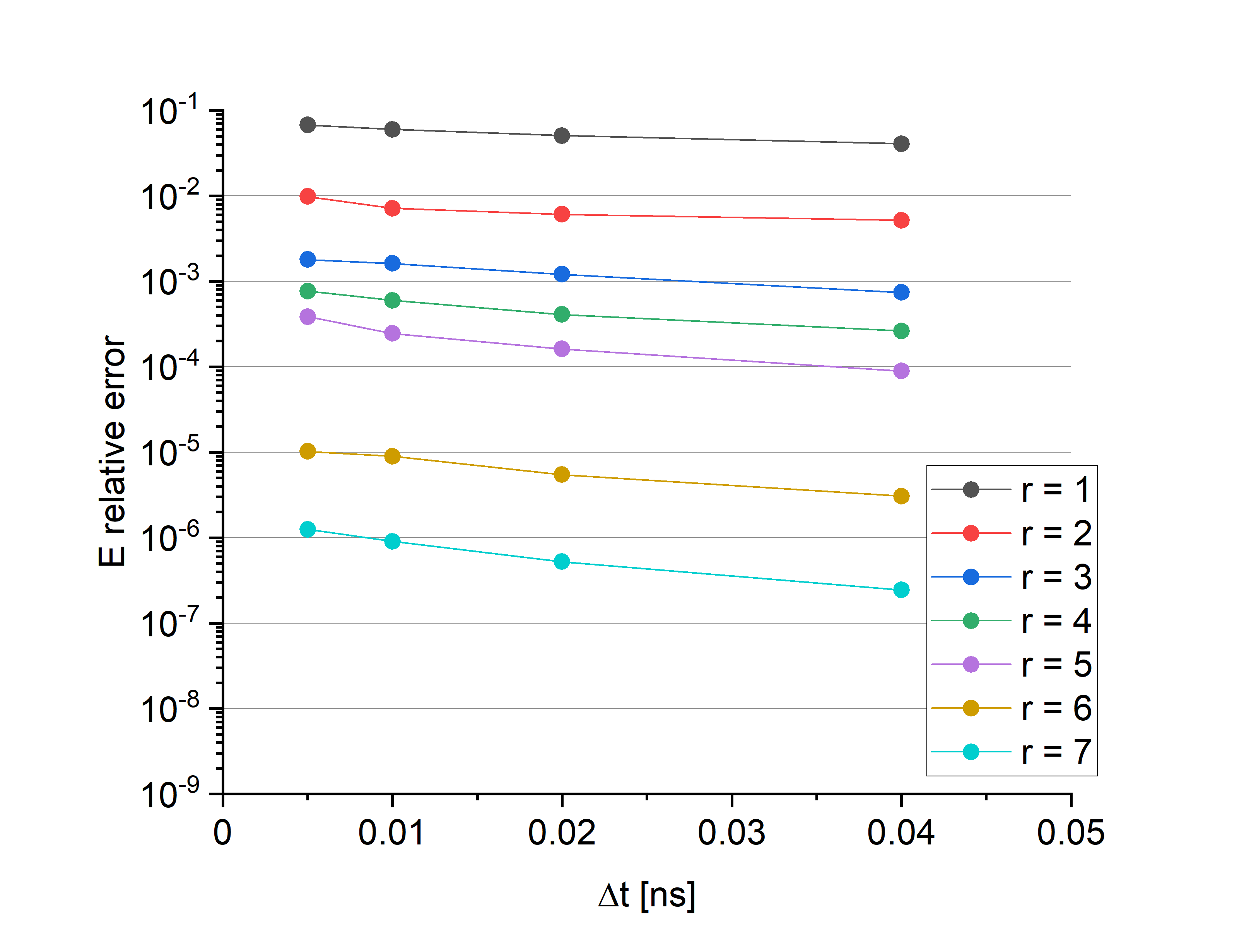}}
	\subfloat[$\frac{||E_h - E_h^r||_{\infty}}{||E_h||_{\infty}}$ at $t=6$ ns]{\includegraphics[width=.35\textwidth]{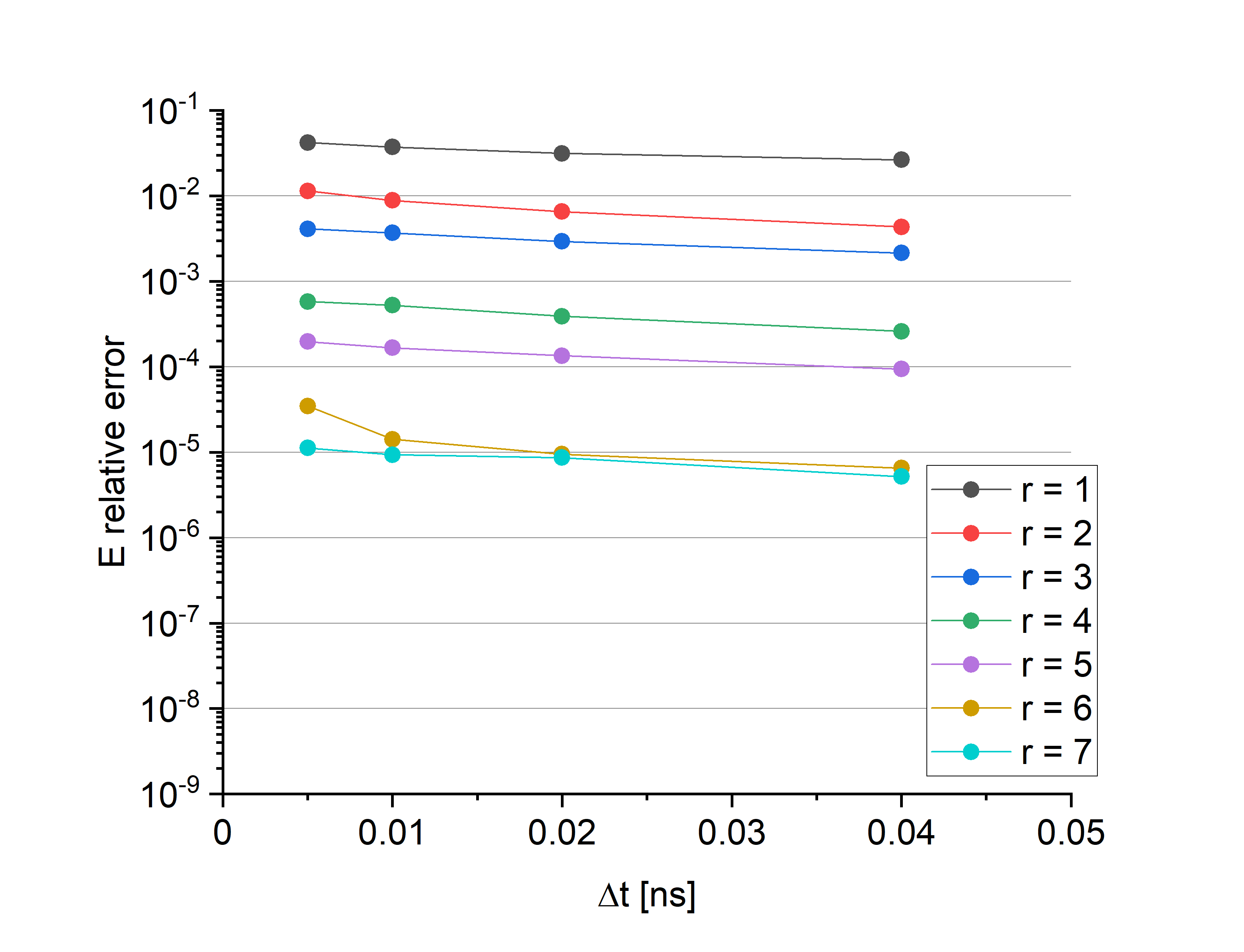}}
	\caption{Results of time step refinement
 the MLQD method with the MBE-SC scheme and POD of the intensity
$\Delta x = 6\times 10^{-2}$ cm. \label{fig:refine-t}}

	\centering \hspace*{-.5cm}
	\subfloat[$\frac{||E_h - E_h^r||_{\infty}}{||E_h||_{\infty}}$ at $t=0.4$ ns]{\includegraphics[width=.35\textwidth]{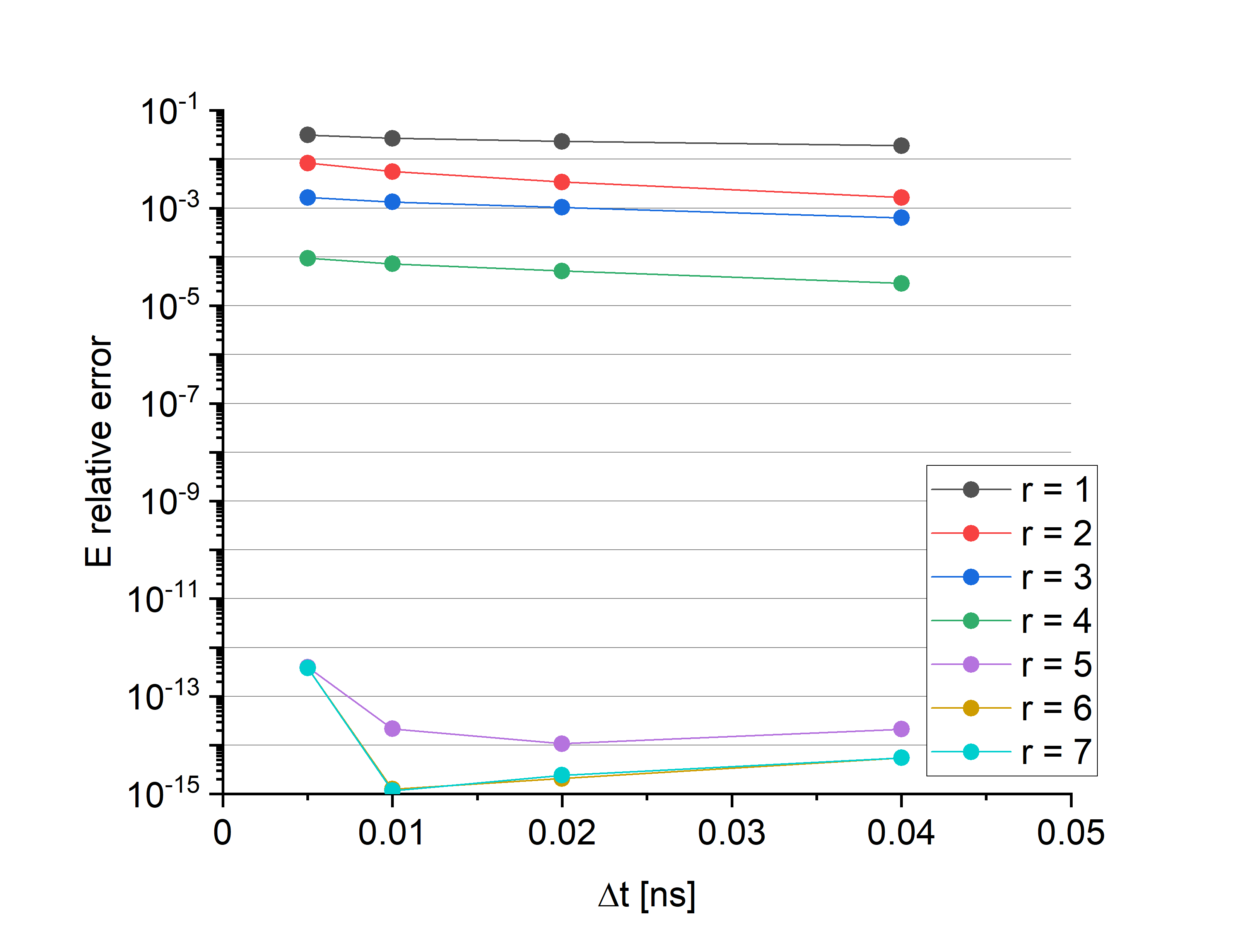}}
	\subfloat[$\frac{||E_h - E_h^r||_{\infty}}{||E_h||_{\infty}}$ at $t=1$ ns]{\includegraphics[width=.35\textwidth]{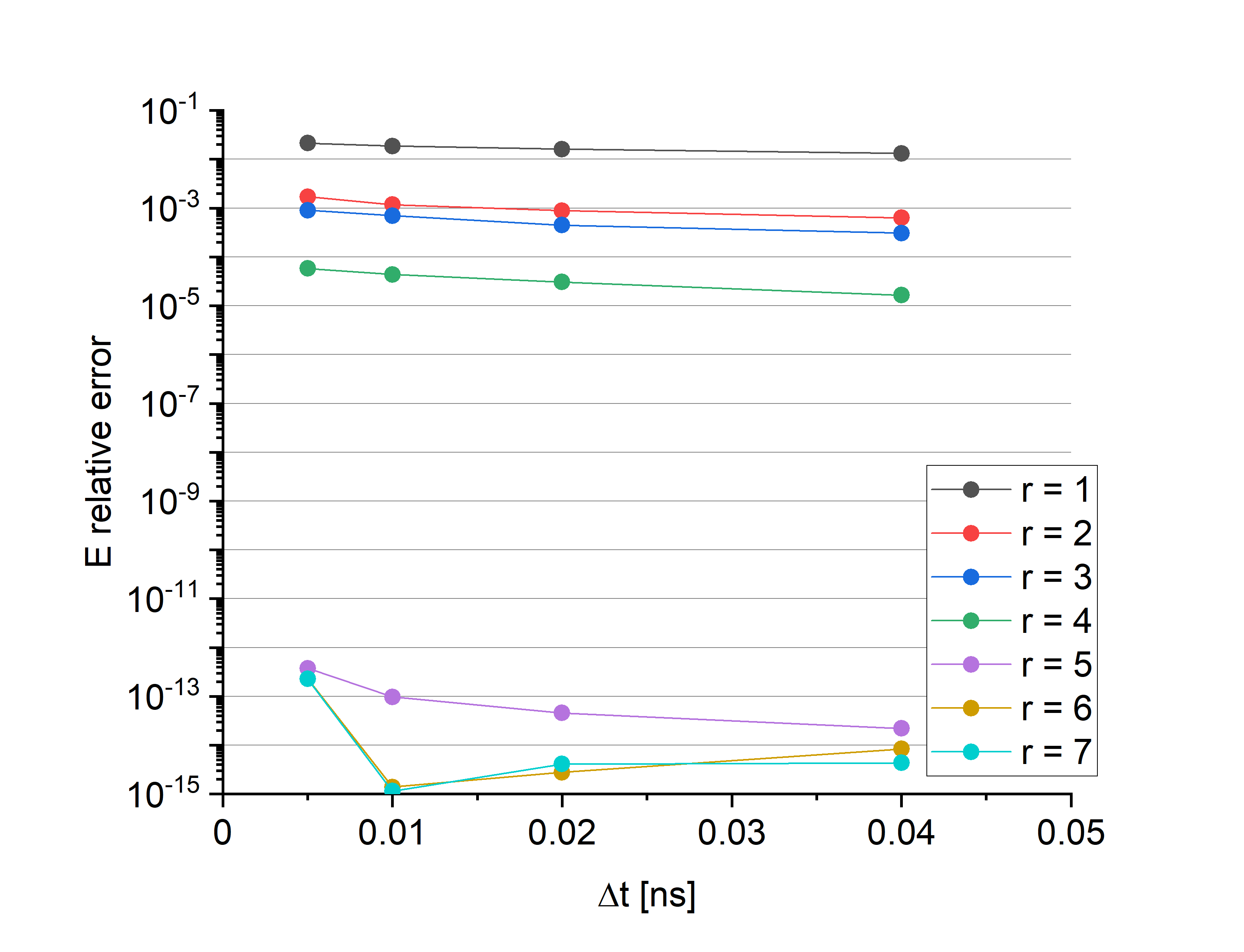}}
	\subfloat[$\frac{||E_h - E_h^r||_{\infty}}{||E_h||_{\infty}}$ at  $t=6$ ns]{\includegraphics[width=.35\textwidth]{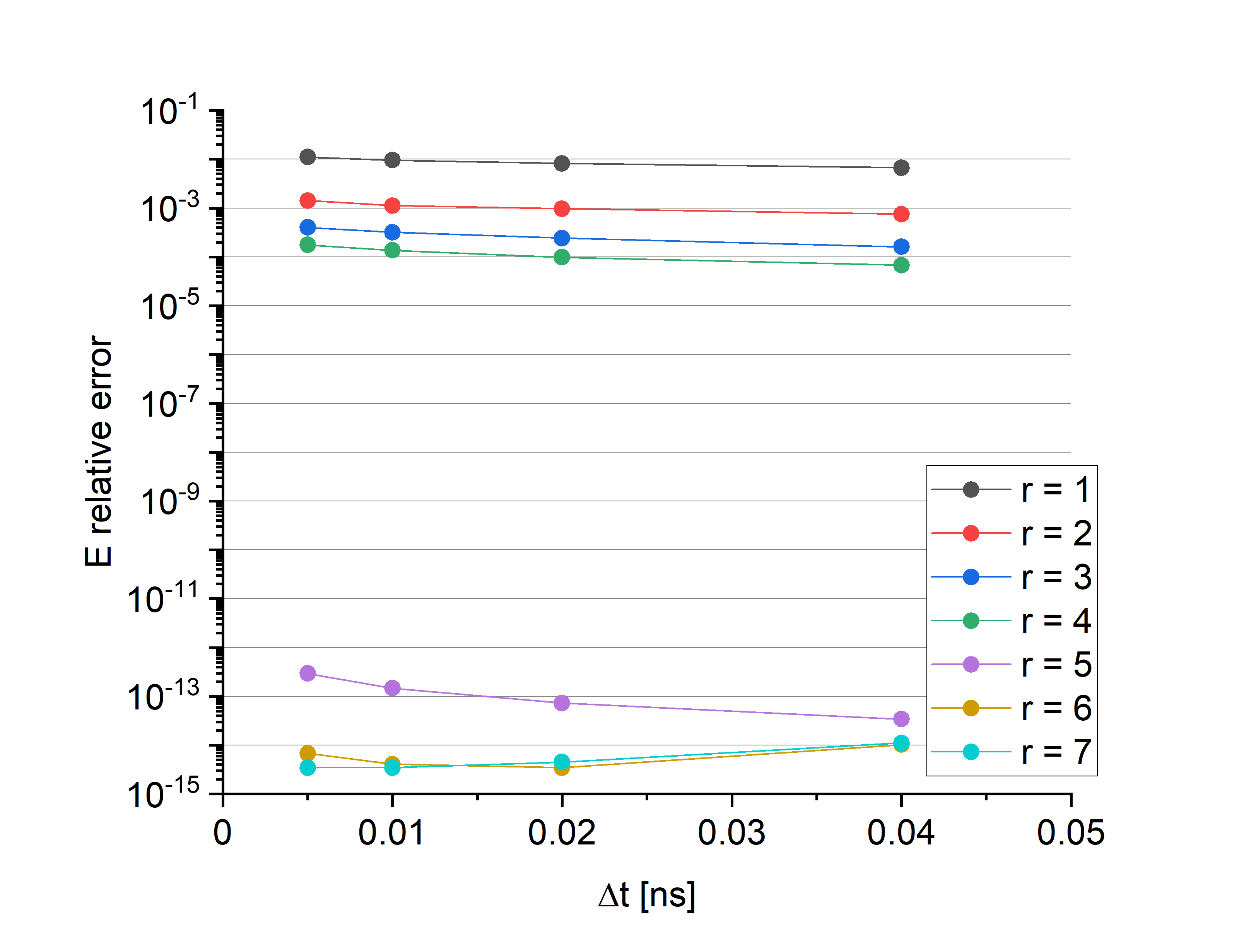}}
	\caption{Results of time step refinement
the MLQD method with the MBE-SC scheme and POD of the remainder term
 for $\Delta x = 6\times 10^{-2}$ cm.
\label{fig:refine-t-P2}}
\end{figure}

 Figures \ref{fig:refine-x} and \ref{fig:refine-x-P2} present the results
  of spatial mesh refinement for the fixed time step size
  \linebreak
  $\Delta t$ = 2$\times $10$^{-2}$ ns.
They show the relative error of E in $\infty$-norm
  for uniform meshes with
    \linebreak
  $\Delta x$ = 0.24, 0.12, 0.06, 0.03  cm.
The number of   degrees of freedom of the discrete intensity increases with  refinement of spatial mesh.
The results show that the change in the relative error decreases with refinement.
The factor of change on fine meshes approaches one.
This indicates that the error due to low-rank POD of data
 representing  intensities tends to a limit as $\Delta x \to 0$ for the fixed time step $\Delta t$.
Figures \ref{fig:refine-t} and \ref{fig:refine-t-P2} present
the relative error of E in $\infty$-norm
for the numerical solution computed with refined  time steps
($\Delta t $ = 4$\times$10$^{-2}$, 2$\times$10$^{-2}$,  10$^{-2}$, 5$\times$10$^{-3}$ ns)
 on the spatial mesh  with $\Delta x$ = 6$\times$10$^{-2}$ cm.
 These results show  increase in the relative error
 in reproducing the discrete solution on the given grids.
More analysis is needed to  study properties of the  methods.

\section{\label{sec:conc} Conclusions}

This paper  presented   implicit methods  with approximate time evolution operator in the high-order
Boltzmann equation and  reduced memory   for TRT problems.
 The obtained  results showed that
 the proposed methods  reproduce the numerical solution of the underlying discretization method on the given phase-space grid
with various degrees of accuracy while  reducing  storage of data between time steps.
The accuracy depends on the rank of the POD of the data representing intensity from the previous time level.
It is possible to achieve accuracy that is good for practical routine simulations and
significantly reduce memory usage.
There are extra computational costs due to calculations of  the POD of intensities.
The proposed approximate implicit methods are intended for  computer architectures on which
one  can take advantage of  extra computations  for reduction of  memory storage.
The proposed approach can be applied to various  time integration methods and different kind of transport problems.

\section*{Acknowledgements}

This research project  is funded  by the Department of Defense, Defense Threat Reduction Agency, grant number HDTRA1-18-1-0042.
The content of the information does not necessarily reflect the position or the policy of the federal government, and no official endorsement should be inferred.

\clearpage

\bibliographystyle{elsarticle-num}
\bibliography{trt-memory-mc2021}

\end{document}